\documentclass[12pt]{amsart} 
\usepackage{amsmath,amssymb}
\usepackage{latexsym}
\usepackage{amssymb} 
\usepackage{latexsym}
\usepackage{delarray}

    \newtheorem{thm}{Theorem}[section]

    \theoremstyle{definition}
    \newtheorem{defn}[thm]{Definition}
    \theoremstyle{remark}

    \numberwithin{equation}{section}

     \newcommand{\Rn}{{{\mathbb R}^n}}

  \newcommand\ve\mathbf

\begin{document}
%
\title[On pseudo-differential operators on group ${\rm SU}(2)$]
    {On pseudo-differential operators on group ${\rm SU}(2)$} 
     \author{Michael Ruzhansky and Ville Turunen}

\address{%
Department of Mathematics\\
Imperial College London\\
United Kingdom}

\email{m.ruzhansky@imperial.ac.uk}

\thanks{The first author was supported by a Royal Society grant and by the EPSRC Grant EP/E062873/01} 

\address{%
Institute of Mathematics\\
Helsinki University of Technology\\
Finland}

\email{ville.turunen@tkk.fi}



\date{December 31, 2007}

\begin{abstract}
In this paper we will outline elements of the global calculus of pseudo-differential operators on the group ${\rm SU}(2)$.
This is a part of a more general approach to pseudo-differential operators on compact Lie groups that will appear in \cite{RT08}.
\end{abstract}

\maketitle

\section{Introduction}

The paper is devoted to outline a global approach to pseudo-differential operators on matrix group ${\rm SU}(2)$, without resorting to local charts.
This can be done by presenting functions on the group by Fourier series obtained from the representations of the group.
Due to non-commutativity,
the Fourier coefficients become matrices of varying dimension.
A pseudo-differential operator can be presented as a convolution operator valued mapping on the group.
The corresponding Fourier coefficient matrices provide a natural global symbol of the pseudo-differential operator.
Consequently, a global calculus and full symbols of pseudo-differential operators can be obtained on the 3-dimensional sphere, with further geometric implications.

Let us first recall and fix the notation for the standard notions of the Fourier analysis.
In view of our applications,
in the formulae below we will emphasise the fact that the dual group $\widehat{\Rn}$ of $\Rn$ is isomorphic to $\Rn$, thus keeping the notation of $\widehat{\Rn}$.
The Fourier transform
$$
   (f\mapsto\widehat f):{\mathcal S}(\Bbb R^n)\to{\mathcal S}(\widehat{\Bbb R^n}) $$ is defined by $$
     \widehat{f}(\xi) = \int_{\Bbb R^n}
           f(x)\ {\rm e}^{-{\rm i}2\pi x\cdot\xi}
           \ {\rm d}x,
$$
and the Fourier inversion formula is
$$
     f(x) = \int_{\widehat{\Bbb R^n}}
           \widehat{f}(\xi)\ {\rm e}^{{\rm i}2\pi x\cdot\xi}
           \ {\rm d}\xi.
$$
Operator $A:{\mathcal S}(\Bbb R^n)\to{\mathcal S}(\Bbb R^n)$ is a pseudo-differential operator of order $m\in\Bbb R$ if $$
     (Af) (x)
           = \int_{\widehat{\Bbb R^n}} \sigma_A(x,\xi)\ \widehat{f}(\xi)
                   \ {\rm e}^{{\rm i}2\pi x\cdot\xi}\ {\rm d}\xi, $$ where the symbol $\sigma_A\in C^\infty(\Bbb R^n\times\widehat{\Bbb R^n})$ of $A$ satisfies symbolic inequalities $$
     \left| \partial_\xi^\alpha\partial_x^\beta\sigma_A(x,\xi) \right|
           \leq C_{A\alpha\beta m}\ \langle\xi\rangle^{m-|\alpha|}, $$ for all multi-indices $\alpha$ and $\beta$ and all $x\in\Rn$ and $\xi\in\widehat{\Rn}$ (see e.g. Kohn and Nirenberg \cite{KN65}, H\"ormander \cite{Ho65}).
In this case we write $A\in\Psi^m(\Bbb R^n)$ and $\sigma_A\in S^m(\Bbb R^n)$.

On a compact manifold $M$ without boundary the space $\Psi^m(M)$ of pseudo-differential operators can be defined via localisations. One disadvantage of this definition is that it often destroys the geometric and other global structures of $M$.
The main problem is that there is no canonical way to write the phase of a pseudo-differential operator when working in local coordinates.
However, if one fixes a connection on the manifold, one can define the phase in a global way in terms of the connection (see e.g. Widom \cite{Wi80}, Safarov \cite{Sa97}, Sharafutdinov \cite{Sh05}).

If $G$ is a Lie group, we may construct another global representation of pseudo-differential operators on $G$ and their calculus by relying on the globally defined group structure.
Moreover, such analysis can be extended to symmetric spaces, where we have a transitive action $G\times M\to M$ of a Lie group $G$ on a manifold $M$.
In this way we obtain the calculus on $M$ as a ``shadow''
from the calculus on $G$.
For example, for $M=\Bbb S^{n-1}$ we can take $G={\rm SO}(n)$, and for the complex sphere $M= \{ x\in\Bbb C^n:\ \|x\|_{\Bbb C^n}=1\}$ we can take $G={\rm SU}(n)$.
A version of pseudo-differential operators on ${\rm SU}(2)$ for functions with finite Fourier series (i.e. for trigonometric polynomials) was analysed in \cite{GK97}.

The important starting point for our analysis is the global Fourier analysis on the Lie group $G$.
Two important ingredients that are needed on $G$ are

(a$_G$) the Fourier transform;

(b$_G$) multiplication and its relation to the Fourier transform.

\medskip
\noindent
Analogies of these operations on $\Bbb R^n$ are

(a$_{\Bbb R^n}$)
$x\mapsto {\rm e}^{{\rm i}x\cdot\xi}\quad (\xi\in\Bbb R^n)$ and the Fourier integral;

(b$_{\Bbb R^n}$) multiplication and the Fourier transform related by $$
     \widehat{(-{\rm i}x)^\alpha f\quad} = \partial_\xi^\alpha \widehat{f}.
$$
On torus $\Bbb T^n = (\Bbb R/2\pi\Bbb Z)^n$ we have the following analogies:

(a$_{\Bbb T^n}$)
$x\mapsto {\rm e}^{{\rm i}x\cdot\xi}\quad (\xi\in\Bbb Z^n)$ and the Fourier series;

(b$_{\Bbb T^n}$) multiplication and the Fourier transform are related by $$
     \widehat{x^{(\alpha)} f} = \triangle_\xi^\alpha \widehat{f}, $$ where for $n=1$, we define $
     x^{(\alpha)} = \left( {\rm e}^{-{\rm i} x} - 1 \right)^\alpha, $ and similarly in higher dimensions.
The relation between the toroidal and Euclidean representations of pseudo-differential operators in one dimension (i.e. on the circle $\Bbb T^1$) was given by Agranovich (\cite{Ag79, Ag85, Ag90}), who showed that $A\in\Psi^m(\Bbb T^1)$ if and only if it can be written as $$
     (Af) (x)
           = \sum_{\xi\in\Bbb Z^1} \sigma_A(x,\xi)\ \widehat{f}(\xi)
                   \ {\rm e}^{{\rm i} x\cdot\xi}, $$ where the toroidal symbol $\sigma_A\in C^\infty(\Bbb T^1\times\Bbb Z^1)$ of $A$ satisfies $$
     \left| \triangle_\xi^\alpha\partial_x^\beta\sigma_A(x,\xi) \right|
           \leq C_{A\alpha\beta m}\ \langle\xi\rangle^{m-|\alpha|}.
$$
Here the difference operator $\triangle_\xi$ is defined by $$
     (\triangle_\xi\sigma)(\xi) = \sigma(\xi+1)-\sigma(\xi).
$$
For $\Psi^m(\Bbb T^n)$ the natural analogy holds with partial difference operators $\triangle_\xi^\alpha$.
This result has been extended to a more comprehensive analysis of pseudo-differential operators and the corresponding toroidal representations of Fourier integral operators in all dimensions in \cite{RT07} and \cite{RT08a}.
We also mention the related paper by Amosov \cite{Am88} on pseudo-differential operators on the circle, and more general characterisations of pseudo-differential operators on the torus by McLean \cite{Mc91}.
We note that  such discrete analysis has applications for numerical mathematics (e.g. \cite{SV02}, \cite{SW87}, \cite{Va99}).

{\bf Acknowledgements:}
The authors wish to express their gratitude to 
Matania~Ben-Artzi for valuable comments leading to an
improvement of the manuscript.

\section{Calculus on compact Lie groups $G$}

A global description and calculus of pseudo-differential operators is possible on general compact Lie groups $G$ and homogeneous spaces of the type $G/K$ with $K\cong\Bbb T^n$.
The main ingredients of this analysis are

\begin{itemize}
\item[(a$_G$)] the Fourier transform based on the usual Haar measure $\mu_G$ on $G$ together with irreducible unitary representations and corresponding Fourier series; \item[(b$_G$)] multiplication and the Fourier series related by the exponential mapping, using Taylor polynomials in deformed exponential coordinates; \item[(c$_G$)] global pseudo-differential operator calculus exploiting the Fourier series.
\end{itemize}
Let us now outline these elements in more detail.

\subsection*{(a$_G$) Fourier transform}

Let $\Bbb C^{c\times d}$ denote the space of complex matrices with $c$ rows and $d$ columns, and let ${\rm U}(d)\subset\Bbb C^{d\times d}$ be the set of $d$-dimensional unitary matrices.
Let ${\mathcal D}(G) = C^\infty(G)$ denote the test function space.
The unitary dual $\widehat G$ consists of equivalence classes $[\xi]$ of irreducible unitary representations $\xi$ of $G$.
In the sequel, from each class $[\xi]\in\widehat G$, we choose one representative $\xi:G\to {\rm U}({\rm dim}(\xi))$, so that the unitary matrices $\xi(x)=\begin{pmatrix} \xi_{mn}(x) \end{pmatrix}_{m,n}$ define ${\rm dim}(\xi)^2$ functions $\xi_{mn}\in C^\infty(G)$.
Fourier coefficient
$\widehat{f}(\xi)\in\Bbb C^{{\rm dim}(\xi)\times {\rm dim}(\xi)}$ of $f\in {\mathcal D}(G)$ is defined by $$
     \widehat{f}(\xi) = \int_G f(x)\ \xi(x)\ {\rm d}\mu_G(x)
     = \begin{pmatrix} \displaystyle
       \int_G f(x)\ \xi_{mn}(x)\ {\rm d}\mu_G(x)
       \end{pmatrix}_{m,n},
$$
so that the Fourier series becomes
\begin{eqnarray*}
     f(x) & = & \sum_{[\xi]\in\widehat{G}} {\rm dim}(\xi)
     \ {\rm Tr}\left( \widehat{f}(\xi)\ \xi^\ast(x) \right), \end{eqnarray*} where $\xi^\ast(x) = \xi(x)^\ast$, and ${\rm Tr}$ is the trace functional (i.e. the sum of the diagonal elements in a matrix).
We can also notice the useful relation
$
     \widehat{f\ast g}(\xi) = \widehat{f}(\xi)\ \widehat{g}(\xi).
$
If $A:{\mathcal D}(G)\to{\mathcal D}(G)$ is a linear continuous map, then it can be written as $$
     Af(x) = \sum_{[\xi]\in\widehat G} {\rm dim}(\xi)
     \ {\rm Tr}\left( \sigma_A(x,\xi)\ \widehat{f}(\xi)\
     \xi^\ast(x) \right),
$$
where the (matrix) symbol of $A$ is defined as the mapping $(x,\xi)\mapsto \sigma_A(x,\xi)\in\Bbb C^{{\rm dim}(\xi)\times{\rm dim}(\xi)}$, where $$
     \sigma_A(x,\xi) = \xi(x)\ (A(\xi^\ast))(x)
     = \begin{pmatrix}\displaystyle
       \sum_{k} \xi_{mk}(x)\ (A\overline{\xi_{nk}})(x)
       \end{pmatrix}_{m,n}.
$$

\subsection*{(b$_G$) multiplication and the Fourier series}

The Taylor expansion of a smooth function $f:G\to\Bbb C$ nearby the neutral element $I\in G$ becomes $$
     f(x) \approx \sum_{|\alpha|\leq N} \frac{1}{\alpha!}\ x^{(\alpha)}
           \ (\partial_x^{(\alpha)} f)(I), $$ with $x^{(\alpha)}\approx x^\alpha$ in e.g. the exponential coordinates.
We can then define the difference operators by the relation $$
    \widehat{x^{(\alpha)}f}(\xi) =: \triangle_\xi^\alpha\widehat{f}(\xi).
$$

\subsection*{(c$_G$) pseudo-differential operators}

The following theorem relates pseudo-differential operators with their globally defined symbols:
\begin{thm}[\cite{Tu01}]\label{TH:symbols}
$A\in{\mathcal L}({\mathcal D}(G))$ belongs to $\Psi^m(G)$ if and only if $$
    \sigma_A\in S^m(G)=\bigcap_{k=0}^\infty S_k^m(G).
$$
\end{thm}
In this theorem, the symbol classes $S^m_k(G)$ are defined in the following recursive way:
\begin{defn}
Symbol $\sigma_A\in S_0^m(G)$ if and only if $$
     \left\| \triangle_\xi^\alpha \partial_x^\beta \sigma_A(x,\xi) \right\|
           \leq C_{A\alpha\beta m}\  \langle \xi \rangle^{m-|\alpha|}, $$ where the norm is the usual operator norm, see Section~\ref{symbolinequalities} for details.
Here, the weight is
$$
   \langle\xi\rangle = \left(1-\lambda_\xi\right)^{1/2}, $$ where $\lambda_\xi\leq 0$ is the eigenvalue of the bi-invariant Laplacian $\Delta:C^\infty(G)\to C^\infty(G)$ on the eigenspace spanned by the matrix element functions $\xi_{mn}\in C^\infty(G)$.
Furthermore, $\sigma_A\in S_{k+1}^m(G)$ if and only if \begin{eqnarray*}
     \sigma_A
     & \in &
           S_k^m(G),\\
     \sigma_{\partial_j}\sigma_A - \sigma_A\sigma_{\partial_j}
     & \in &
           S_k^m(G),\\
     (\triangle_\xi^\gamma\sigma_{\partial_j})\ \sigma_A
     & \in &
           S_k^{m+1-|\gamma|}(G),\\
     (\triangle_\xi^\gamma\sigma_A)\ \sigma_{\partial_j}
     & \in &
           S_k^{m+1-|\gamma|}(G),
\end{eqnarray*}
where $|\gamma|>0$ and $1\leq j\leq{\rm dim}(G)$.
Then we define
$$
     S^m(G) = \bigcap_{k=0}^\infty S_k^m(G).
$$
\end{defn}
We note that the proof of Theorem \ref{TH:symbols} is based on the commutator characterisation of pseudo-differential operators (see e.g. \cite{Be77}, \cite{Du77}, \cite{Co79}, \cite{Tu00}; see also \cite{Ta97} for a different characterisation on the sphere, studying smoothness of operator valued mappings).
It is well-known that if $A\in\Psi^{m_A}$ and $B\in\Psi^{m_B}$ then $AB,BA\in \Psi^{m_A+m_B}$, but the commutator satisfies $$
          \quad [A,B]=AB-BA\in\Psi^{m_A+m_B-1}.
$$
The commutator characterisation uses
this property to characterise pseudo-differential operators:
\begin{thm}
$A\in\Psi^m(M)$ if and only if
$$
     [D_{k+1},[D_k,\cdots[D_1,A]\cdots]]\in {\mathcal L}(H^m(M),H^0(M)) $$ for every sequence of smooth vector fields $D_k$ on $M$, where $H^m(M)$ is the Sobolev space of order $m\in\Bbb R$ on $M$.
\end{thm}
Details of these constructions with applications to the analysis on Lie groups will appear in \cite{RT08}.
We also note that this applies also to non-invariant pseudo-differential operators (compared to e.g. \cite{St71}, \cite{St72}).

\section{Calculus on $G={\rm SU}(2)$}

Let us now explain the constructions of the previous section in more detail in the case of the matrix group ${\rm SU}(2)$.
In this case the special unitary group of complex rotations $$
     G = {\rm SU}(2)
     = \{x\in\Bbb C^{2\times 2}:\ x^\ast = x^{-1},\ {\rm det}(x)=1\} $$ acts transitively on the unit sphere in $\Bbb C^2$, so the analysis also applies for pseudo-differential operator on spheres in $\Bbb C^2$.
The group ${\rm SU}(2)$ is especially important since $G\cong \Bbb S^3$, the 3-sphere with the quaternionic structure.
In particular, we obtain the global calculus of pseudo-differential operators on $\Bbb S^3$, which allows further extensions to general closed simply connected 3-manifolds by the pullback by the Ricci flow (details of this will appear in \cite{RT08b}).

Let us define basic rotations
$\omega_j(t)\in G$ ($j=1,2,3$) by
$$
           \begin{pmatrix}
                   \cos\frac{t}{2} & {\rm i}\sin\frac{t}{2} \\
                   {\rm i}\sin\frac{t}{2} & \cos\frac{t}{2}
           \end{pmatrix},\quad
           \begin{pmatrix}
                   \cos\frac{t}{2} & -\sin\frac{t}{2} \\
                   \sin\frac{t}{2} & \cos\frac{t}{2}
           \end{pmatrix},\quad
           \begin{pmatrix}
                   {\rm e}^{{\rm i}t/2} & 0 \\
                   0 & {\rm e}^{-{\rm i}t/2}
           \end{pmatrix}.
$$
The Euler angles $\phi,\theta,\psi$
of $x\in G$ satisfy
$$
     x = x(\phi,\theta,\psi)
           := \omega_3(\phi)\ \omega_2(\theta)\ \omega_3(\psi), $$ where $-\pi<\phi\leq \pi$, $0\leq\theta\leq\pi$, $-2\pi<\psi\leq 2\pi$.
The Euler angles define local coordinates for $G$ whenever $0<\theta<\pi$ (with natural interpretation when $\phi=\pi$ or $\psi=2\pi$).
Let us now define the principal ingredients for the analysis.

\subsection*{(a$_G$) Fourier transform}

The Haar integral on $G={\rm SU}(2)$ in the Euler angles is given by \begin{eqnarray*}
     f
           & \mapsto &
           \int_G f\ {\rm d}\mu_G\\
           & = &
           \frac{1}{16\pi^2}\int_{-2\pi}^{2\pi} \int_0^\pi \int_{-\pi}^{\pi}
           f(x)
           \ \sin(\theta)\ {\rm d}\phi\ {\rm d}\theta\ {\rm d}\psi.
\end{eqnarray*}
The unitary dual of $G={\rm SU}(2)$ is
$$
     \widehat{G}\cong \left\{t^l\in {\rm Hom}(G,{\rm U}(2l+1)):
     \ l\in\frac{1}{2}\Bbb N \right\},
$$
where
${\rm U}(d)\subset\Bbb C^{d\times d}$
is the unitary matrix group of dimension $d$, and functions $t^l_{mn}\in C^\infty(G)$ are products of exponentials and Legendre functions.

The Fourier transform $f\mapsto\widehat{f}$ on $G={\rm SU}(2)$ becomes $$
     \widehat{f}(l) = \widehat{f}(t^l) := \int_G f(x)\ t^l(x)\ {\rm d}\mu_G(x), $$ with the inverse Fourier transform given by $$
     f = \sum_{l\in\Bbb N/2} (2l+1)
     \ {\rm Tr}\left(\widehat{f}(l)\ (t^l)^\ast \right).
$$
For details of these formulae see
e.g. Vilenkin \cite{Vi68} and Zelobenko \cite{Ze73}.
We also note that the Fourier transform
$\widehat{f}(\xi)$ of $f\in{\mathcal D}(G)$ is defined for every $\xi\in[\xi]\in\widehat{G}$, i.e.
$\widehat{f}(\xi):V_\xi\to V_\xi$
is linear on the representation vector space $V_\xi$ of the irreducible representation $\xi:G\to{\rm Aut}(V_\xi)$ of $G$.
If $\xi\sim\eta$ (i.e. $[\xi]=[\eta]\in\widehat G$) then $\widehat{f}(\xi)$ and $\widehat{f}(\eta)$ are related by intertwining.

\subsection*{(b$_G$) Taylor polynomials}

Lie group $G={\rm SU}(2)$ has
Lie algebra ${\frak g}={\rm su}(2)=\{X(z):\ z\in\Bbb R^3\}$, where $$
     X(z) = \frac{1}{2} \begin{pmatrix}
           {\rm i} z_3   & {\rm i} z_1 - z_2 \\
           {\rm i} z_1 + z_2&  -{\rm i} z_3 \\
           \end{pmatrix}.
$$
If $\{e_j\}_{j=1}^3$ is the standard basis of $\Bbb R^3$ and $X_j=X(e_j)$ (Pauli matrices) then we have commutator relations:
$$
     [X_1,X_2]=X_3,\quad
     [X_2,X_3]=X_1,\quad
     [X_3,X_1]=X_2.
$$
Now, let us describe the vector fields in the basis of the Lie algebra.
For $Y\in{\rm su}(2)$, let
$$
     A_Y:C^\infty({\rm SU}(2))\to C^\infty({\rm SU}(2)) $$ be the left-invariant differential operator defined by \begin{equation*}
     A_Y f (x)
     = \left. \frac{\rm d}{{\rm d}t} f(x\ \exp(tY)) \right|_{t=0}.
\end{equation*}
Let us denote $A_j := A_{X_j}$.
For a multi-index $\beta\in\Bbb N^3$,
define the left-invariant differential operator $$
   \partial^\beta := A_1^{\beta_1} A_2^{\beta_2} A_3^{\beta_3}.
$$
Let us also define creation, annihilation, and neutral operators by the following formulae:
\begin{eqnarray*}
     \partial_+ & := & {\rm i} A_1 - A_2, \\
     \partial_- & := & {\rm i} A_1 + A_2, \\
     \partial_0 & := & {\rm i} A_3.
\end{eqnarray*}
In Euler angles these operators can be expressed as \begin{eqnarray*}
     \partial_+ & = & {\rm e}^{-{\rm i}\psi}
     \left( {\rm i}\ \frac{\partial}{\partial\theta}
           -\frac{1}{\sin(\theta)}\ \frac{\partial}{\partial\phi}
           +\frac{\cos(\theta)}{\sin(\theta)}
           \ \frac{\partial}{\partial\psi} \right), \\
     \partial_- & = & {\rm e}^{{\rm i}\psi}
     \left( {\rm i}\ \frac{\partial}{\partial\theta}
           +\frac{1}{\sin(\theta)}\ \frac{\partial}{\partial\phi}
           -\frac{\cos(\theta)}{\sin(\theta)}
           \ \frac{\partial}{\partial\psi} \right), \\
     \partial_0 & = & {\rm i} \frac{\partial}{\partial\psi}.
\end{eqnarray*}
Note that the invariant Laplace operator on $G$ is given by $$
     \Delta = A_1^2 + A_2^2 + A_3^2.
$$

Let $t:=\|z\|_{\Bbb R^3}/2$. The Rodrigues formula states that we have $$
     \exp(X(z)) = \begin{pmatrix} 1 & 0 \\ 0 & 1\end{pmatrix} \cos(t)
     + X(z) \frac{\sin(t)}{t}.
$$
Define now functions $x_{ij}\in C^\infty(G)$ by $$
     x=\begin{pmatrix} x_{11}(x) & x_{12}(x) \\ x_{21}(x) & x_{22}(x)
     \end{pmatrix} \in G.
$$
These functions (or their linear combinations) will play the role of coordinates.
We can now define the Taylor polynomial
$x^{(\alpha)}\in C^\infty(G)$ ($\alpha\in\Bbb N^3$) by setting \begin{eqnarray*}
     x^{(\alpha)} & = &
     \left(x_{12}\right)^{\alpha_1}\ \left(x_{21}\right)^{\alpha_2}
     \ \left(x_{11}-x_{22}\right)^{\alpha_3}.
\end{eqnarray*}
The corresponding difference operators $\triangle_\xi^\alpha$ are then defined by $$
     \triangle_\xi^\alpha \widehat{f}(\xi)
     := \widehat{x^{(\alpha)}f}(\xi).
$$
In the Taylor expansion formula, Taylor polynomials are multiplied by the partial differential operators of corresponding orders.
For this purpose, we define
$$
     \partial^{(\alpha)} = \sum_{|\beta|\leq |\alpha|} c_{\alpha\beta}
           \ \partial^\beta,
$$
so that we have
$$
     \partial^{(\alpha)} x^{(\beta)}|_{x=I} = \alpha!\
     \delta_{\alpha,\beta}.
$$
Now the Taylor expansion for a function $f\in{\mathcal D}(G)$ at the neutral element $I\in G$ becomes $$
     f(x) \sim \sum_{\alpha\geq 0} \frac{1}{\alpha!}\ x^{(\alpha)}
           \ (\partial^{(\alpha)}f)(I).
$$
Let us finally briefly recall
the irreducible unitary representations for the group ${\rm SU}(2)$.
Group ${\rm SU}(2)$ acts naturally on space $\Bbb C[z_1,z_2]$ by the formula $$
     T:{\rm SU}(2)\to {\rm GL}(\mathbb C[z_1,z_2]),\quad
     (T(x)f)(z) := f(zx).
$$
Invariant subspaces $V_l\subset\Bbb C[z_1,z_2]$ of homogeneous $2l$-degree polynomials, where $l\in\Bbb N/2$, give the restriction $$
     T_l:{\rm SU}(2)\to {\rm GL}(V_l),\quad (T_l(x) f)(z) = f(zx).
$$
If $T_\infty$ is an irreducible representation of ${\rm SU}(2)$ on a vector space $V_\infty$ then $T_\infty$ is equivalent to $T_l$, where ${\rm dim}(V_\infty) = 2l+1$.
Operators $T_l(x)$ can be viewed as matrices and now we will recall the matrix elements of $T_l$, which can be written as $$
     t^l_{mn}(\phi,\theta,\psi) =
     {\rm e}^{-{\rm i}(m\phi+n\phi)}\ P^l_{mn}(\cos(\theta)), $$ where $$
     P^l_{mn}(x)  =  c^l_{mn}
     \frac{(1-x)^{(n-m)/2}}{(1+x)^{(m+n)/2}}
      \left(\frac{{\rm d}}{{\rm d}x}\right)^{l-m}
     \left[ (1-x)^{l-n} (1+x)^{l+n} \right], $$ with constants $$
     c^l_{mn} = 2^{-l} \frac{(-1)^{l-n}\ {\rm i}^{n-m}}{\sqrt{(l-n)!\ (l+n)!}}
     \ \sqrt{\frac{(l+m)!}{(l-m)!}}.
$$
Multiplication of two elements $t^{l'}_{m'n'}$ and $t^l_{mn}$ can be expressed as $$
    t^{l'}_{m'n'}\ t^l_{mn} =
     \sum_{k=|l-l'|}^{l+l'} C^{l l' (l+k)}_{m' m (m'+m)}
           \ C^{l l' (l+k)}_{n' n (n'+n)}
           \ t^{l+k}_{(m'+m)(n'+n)},
$$
with the Clebsch--Gordan coefficients $C^{l l' (l+k)}_{m' m (m'+m)}$.
In particular, let us denote
$$
    \begin{pmatrix} t_{--} & t_{-+} \\ t_{+-} & t_{++} \end{pmatrix}
    = t^{1/2} = \begin{pmatrix}
     t^{1/2}_{-1/2,-1/2} & t^{1/2}_{-1/2,+1/2} \\
     t^{1/2}_{+1/2,-1/2} & t^{1/2}_{+1/2,+1/2}
    \end{pmatrix}
$$
and $x^{\pm} := x\pm 1/2$ for $x\in\mathbb R$.
Then we have the following multiplication formulae \begin{eqnarray*}
(2l+1)  t^l_{mn} t_{--}
     & = &
       \sqrt{(l-m+1)(l-n+1)}\ t^{l^+}_{m^- n^-}
     - \sqrt{(l+m)(l+n)}\ t^{l^-}_{m^- n^-}, \\
(2l+1)  t^l_{mn} t_{++}
     & = &
       \sqrt{(l+m+1)(l+n+1)}\ t^{l^+}_{m^+ n^+}
     - \sqrt{(l-m)(l-n)}\ t^{l^-}_{m^+ n^+}, \\
(2l+1)  t^l_{mn} t_{-+}
     & = &
       \sqrt{(l-m+1)(l+n+1)}\ t^{l^+}_{m^- n^+}
     - \sqrt{(l+m)(l-n)}\ t^{l^-}_{m^- n^+}, \\
(2l+1)  t^l_{mn} t_{+-}
     & = &
       \sqrt{(l+m+1)(l-n+1)}\ t^{l^+}_{m^+ n^-}
     - \sqrt{(l-m)(l+n)}\ t^{l^-}_{m^+ n^-}.
\end{eqnarray*}

\section{Differences for symbols on ${\rm SU}(2)$}

Let us define functions
\begin{eqnarray}
     q_+ &:=& t_{-+}, \\
     q_- &:=& t_{+-}, \\
     q_0 &:=& t_{++}-t_{--}.
\end{eqnarray}
The key property of these functions is that each $q_+,q_-,q_0\in {\mathcal D}(G)$ vanishes at the neutral element $I\in G$.

Let $\sigma=\widehat{s}$, which means that $\sigma:{\mathcal D}(G)\to {\mathcal D}(G)$ is the left convolution operator with convolution kernel $s\in {\mathcal D}'(G)$:
$$
     \sigma f = s\ast f.
$$
Let us define ``difference operators'' $\triangle_+,\triangle_-,\triangle_0$
acting on symbols $\sigma$ by
\begin{eqnarray}
     \triangle_+ \sigma &:=& \widehat{q_+\ s}, \\
     \triangle_- \sigma &:=& \widehat{q_-\ s}, \\
     \triangle_0 \sigma &:=& \widehat{q_0\ s}.
\end{eqnarray}
These difference operators will play the role of derivatives in the definition and calculus of pseudo-differential operators.
In general, they depend on the structure of the dual space $\widehat G$ of the group $G$.
For example, in the Euclidean space we have $\widehat{\Bbb R^n}\cong\Bbb R^n$, so their analogues are simply the differential operators $\partial_\xi^\alpha$.
For the torus, we have $\widehat{\Bbb T^n}\cong\Bbb Z^n$, and the corresponding difference operators $\triangle_\xi^\alpha$ have been analysed e.g. in \cite{RT07} and \cite{Tu00}.

We can now define the conditions on
globally defined pseudo-differential symbols on ${\rm SU}(2)$ in terms of these difference operators.
First let us show that a pseudo-differential operator can be globally parametrised by a matrix symbol $\sigma_A(x,\xi)$, where $x\in{\rm SU}(2)$, $\xi$ is a discrete variable corresponding to the dual of ${\rm SU}(2)$, and the size of the matrix $\sigma_A(x,\xi)$ increases with $\xi$.
In fact, using the structure of the dual of ${\rm SU}(2)$, we can identify $\xi$ with a half-integer in $\Bbb N/2$, so that $\sigma_A(x,\xi)$ becomes a matrix of the size $(2\xi+1)\times(2\xi+1)$.
In general, the matrix symbol of a continuous linear operator $A:{\mathcal D}(G)\to{\mathcal D}(G)$ will be written as $\sigma_A(x,\xi)$, and it is defined for every $x\in G$ and $\xi\in[\xi]\in\widehat G$.
If $[\xi]=[\eta]$ then $\sigma_A(x,\xi)$ and $\sigma_A(x,\eta)$ are related and allow a geometric interpretation.
We refer to \cite{RT08b} for further details, and will only give here the description based on the identification of the dual of ${\rm SU}(2)$ with the set $\Bbb N/2$ of half-integers.

By using the Fourier series,
we can conclude that a continuous linear operator $A:{\mathcal D}(G)\to{\mathcal D}(G)$ belongs to $\Psi^m(G)$ if and only if it can be written as \begin{eqnarray*}
     (Af)(x)
     & = & \sum_{\xi\in\Bbb N/2} (2\xi+1)\  {\rm Tr}\left(
                   \sigma_A(x,\xi)\ \widehat{f}(\xi)\
                   t^\xi(x)^\ast \right) \\
     & = & \sum_{\xi\in\Bbb N/2} (2\xi+1) \sum_{m=-\xi}^\xi\sum_{n=-\xi}^\xi
                   \overline{t^\xi_{mn}(x)}
            \left( \sum_{k=-\xi}^\xi \sigma_A(x,\xi)_{mk}
                   \ \widehat{f}(\xi)_{kn}\right), \end{eqnarray*} where the symbol $$
    \sigma_A(x,\xi) = t^\xi(x)\ (A(t^\xi)^\ast)(x)
    = \begin{pmatrix}\displaystyle
      \sum_{k=-\xi}^{\xi} t^\xi_{mk}(x)\ (A\overline{t^\xi_{nk}})(x)
      \end{pmatrix}_{m,n}
      = \begin{pmatrix} \sigma_A(x,\xi)_{mn} \end{pmatrix}_{m,n} $$ satisfies symbol inequalities that will be described below.
With a natural embedding
${\mathcal D}(\Bbb S^2)\hookrightarrow {\mathcal D}(G)$, we also note that every pseudo-differential operator $B\in \Psi^m(\Bbb S^2)$ has an extension $A\in\Psi^m(G)$ such that $A|_{{\mathcal D}(\Bbb S^2)} = B$, see \cite{Tu01}.

Let us consider the following example:
for $\xi\in\Bbb N/2$ and a $(2\xi+1)\times (2\xi+1)$ matrix $A$, let $A_{mn}$ denote the matrix element on the $m$th row and $n$th column, where $|m|,|n|\leq \xi$ such that $m,n\in\{\pm \xi,\pm (\xi-1),\pm (\xi-2),\cdots\}$.
The symbols of the first order partial differential operators $\partial_+,\partial_-,\partial_0$ are then given by $$
     \sigma_{\partial_+}(x,\xi)_{mn} =
     \begin{cases}
       -\sqrt{(\xi-n)(\xi+n+1)}, & {\rm if}\ m=n+ 1, \\
       0, & {\rm otherwise}.
       \end{cases}
$$
$$
     \sigma_{\partial_-}(x,\xi)_{mn} =
     \begin{cases}
       -\sqrt{(\xi+n)(\xi-n+1)}, & {\rm if}\ m=n- 1, \\
       0, & {\rm otherwise}.
       \end{cases}
$$
$$
     \sigma_{\partial_0}(x,\xi)_{mn} =
     \begin{cases}
       n, & {\rm if}\ m=n, \\
       0, & {\rm otherwise},
       \end{cases}
$$
$$
     \sigma_{I}(x,\xi)_{mn} =
     \begin{cases}
       1, & {\rm if}\ m=n, \\
       0, & {\rm otherwise},
       \end{cases}
$$
where $I=(f\mapsto f):{\mathcal D}(G)\to{\mathcal D}(G)$ is the identity operator.
Moreover, we have the following properties:
$$
     \triangle_+ \sigma_{\Delta} = -\sigma_{\partial_-},\quad
     \triangle_- \sigma_{\Delta} = -\sigma_{\partial_+},\quad
     \triangle_0 \sigma_{\Delta} = -\sigma_{\partial_0}, $$ where $\Delta$ is the invariant Laplacian of $G$, $$
     \sigma_I =
     \triangle_+ \sigma_{\partial_+} =
     \triangle_- \sigma_{\partial_-} =
     \triangle_0 \sigma_{\partial_0},
$$
and
$$
     0 =
     \triangle_+ \sigma_I =
     \triangle_- \sigma_I =
     \triangle_0 \sigma_I.
$$
Furthermore,
$
     \triangle_\alpha \sigma_{\partial_\eta}(x,\xi) = 0, $ if $\alpha,\eta\in\{+,-,0\}$ are such that $\alpha\not=\eta$.

\section{Symbol inequalities}\label{symbolinequalities}

In this section we study inequalities describing symbols of pseudo-differential operators on $G={\rm SU}(2)$.
For a vector $v=(v_j)_{j=1}^n\in\Bbb C^d$ we use the Euclidean norm $\|v\|_{\Bbb C^d}$ given by $$
    \|v\|_{\Bbb C^d}^2 := \sum_{j=1}^d |v_j|^2, $$ and for a matrix $M\in\Bbb C^{d\times d}$ the corresponding operator norm $$
    \|M\| = \|M\|_{\Bbb C^{d\times d}}
    := \sup \left\{ \|Mv\|_{\Bbb C^d}:\ v\in\Bbb C^d,\ \|v\|_{\Bbb C^d}\leq 1
      \right\}.
$$
We notice that we have the operator norm equality $$
    \left\| f\mapsto f\ast a \right\|_{{\mathcal L}(L^2(G))}
    = \sup \left\{ \left\|\widehat{a}(\xi)\right\|:\ \xi\in\widehat{G} \right\}.
$$
In order to describe matrix symbols of
standard pseudo-differential operators on ${\rm SU}(2)$, we need to check conditions of Theorem \ref{TH:symbols}.
First, $\sigma_A \in S_0^m(G)$, i.e. it satisfies \begin{equation}\label{EQ:orderm}
    \left\| \triangle_\xi^\alpha \partial_x^\beta \sigma_A(x,\xi) \right\|
     \leq C_{\alpha\beta}\ \langle\xi\rangle^{m-|\alpha|} \end{equation} for all multi-indices $\alpha, \beta$, all $x\in{\rm SU}(2)$, and all $\xi\in\Bbb N/2$.
On $G={\rm SU}(2)$, the weight $\langle\xi\rangle$ becomes $$
   \langle\xi\rangle = \left(1 + \xi + \xi^2\right)^{1/2}.
$$
We note that condition
$(\triangle_\xi^\gamma\sigma_{\partial_j}) \sigma_A \in S_k^{m+1-|\gamma|}(G)$ is automatically satisfied due to the example presented at the end of the previous section.

We will now show that in order to
satisfy conditions of Theorem \ref{TH:symbols}, matrices $\sigma_A(x,\xi)$ must have a certain ``rapid off-diagonal decay'' property.
By the example in the previous section,
$\sigma_{\partial_0}(x,\xi)_{ij} = i\ \delta_{ij}$, so that \begin{eqnarray*}
    && [\sigma_{\partial_0}(x,\xi),\sigma_A(x,\xi)]_{ij} \\
    & = &  \sum_k \left(
      \sigma_{\partial_0}(x,\xi)_{ik}\ \sigma_A(x,\xi)_{kj}
      - \sigma_A(x,\xi)_{ik}\ \sigma_{\partial_0}(x,\xi)_{kj} \right) \\
    & = &  \sum_k \left( i\ \delta_{ik}\ \sigma_A(x,\xi)_{kj}
      - \sigma_A(x,\xi)_{ik}\ k\ \delta_{kj} \right) \\
    & = & (i-j)\ \sigma_A(x,\xi)_{ij}.
\end{eqnarray*}
Let us iterate such a commutator $p\in\Bbb Z^+$ times, to obtain that the symbol $$
    (x,\xi) \mapsto \begin{pmatrix}
      (i-j)^p \sigma_A(x,\xi)_{ij} \end{pmatrix}_{i,j} $$ must belong to $S^m(G)\subset S^m_0(G)$, regardless of $p$.
Thus the pseudo-differential symbol $\sigma_A$ must have ``{\it rapid off-diagonal decay}'':
$$
    \sup_{\xi,i,j}\ \langle\xi\rangle^{-m} \langle i-j\rangle^p
    \left|\sigma_A(x,\xi)_{ij}\right| < \infty $$
--- notice that in a matrix
$\begin{pmatrix} a_{ij} \end{pmatrix}_{i,j}$, the distance from the location of the matrix element $a_{ij}$ to the diagonal is $|i-j|$.

What about the commutators
$[\sigma_{\partial_+},\sigma_A]$ and $[\sigma_{\partial_-},\sigma_A]$?
Due to the symmetries of the symbols of $\partial_\pm$, studying the case of $\sigma_{\partial_+}$ will be enough.
By the example in the previous section,
$$
    \sigma_{\partial_+}(x,\xi)_{ij} = -\sqrt{(\xi-i)(\xi+i+1)}\ \delta_{i+1,j}, $$ so that \begin{eqnarray*}
    && [\sigma_{\partial_+}(x,\xi),\sigma_A(x,\xi)]_{ij} \\
    & = &  \sum_k \left(
      \sigma_{\partial_+}(x,\xi)_{ik}\ \sigma_A(x,\xi)_{kj}
      - \sigma_A(x,\xi)_{ik}\ \sigma_{\partial_+}(x,\xi)_{kj} \right) \\
    & = &  \sum_k \left( -\sqrt{(\xi-i)(\xi+i+1)}\ \delta_{i+1,k}
      \ \sigma_A(x,\xi)_{kj} \right. \\
      && \quad\quad \left.
      + \sigma_A(x,\xi)_{ik}\ \sqrt{(\xi-k)(\xi+k+1)}\ \delta_{k+1,j}
\right) \\
    & = & -\sqrt{(\xi-m)(\xi+m+1)} \ \sigma_A(x,\xi)_{i+1,j}  \\
      &&
      + \sqrt{(\xi-n+1)(\xi+n)}\ \sigma_A(x,\xi)_{i,j-1}.
\end{eqnarray*}
At first sight, this commutator may look baffling:
is this some sort of weighted difference operator acting on $\sigma_A$ ``along the diagonal''?
We have to understand this first-order commutator condition properly before we may consider higher order iterated commutators.

So how to understand the behaviour of
$[\sigma_{\partial_+}(x,\xi),\sigma_A(x,\xi)]$?
We claim that this is a variant of the
``rapid off-diagonal decay'' property.
Let us explain what we mean by this.
Briefly: $\partial_0={\rm i} A_3$,
so the $A_3$-symbol commutator condition means ``rapid off-diagonal decay''.
``Badly behaving'' operators $\partial_+,\partial_-$ are linear combinations of left-invariant vector fields $A_1,A_2$, which are conjugates to $A_3$ (and hence essentially similar in behaviour); then the idea is to exploit the diffeomorphism invariance of $\Psi^m(G)$.

We know from the local theory that the pseudo-differential operator class $\Psi^m(G)$ is diffeomorphism-invariant:
if $\phi:G\to G$ is a diffeomorphism and $A\in\Psi^m(G)$, then also $A_\phi\in\Psi^m(G)$, where $$
    A_\phi f = A(f\circ\phi)\circ\phi^{-1}.
$$
Especially, let us consider the inner automorphisms $$
    \phi_u=(x\mapsto u^{-1} x u):G\to G, $$ where $u\in G$.
Such $\phi_u:G\to G$ is a diffeomorphism that maps one-parametric subgroups to one-parametric subgroups.
If $\xi\in[\xi]\in\widehat{G}$ then
\begin{eqnarray*}
    \widehat{f\circ\phi_u}(\xi)
    & = & \int_G f(u^{-1}xu)\ \xi(x)\ {\rm d}\mu_G(x) \\
    & = & \int_G f(x)\ \xi(uxu^{-1})\ {\rm d}\mu_G(x) \\
    & = & \xi(u) \int_G f(x)\ \xi(x)\ {\rm d}\mu_G(x)\ \xi(u)^\ast \\
    & = & \xi(u)\ \widehat{f}(\xi)\ \xi(u)^\ast, \end{eqnarray*} which is similar to $\widehat{f}(\xi)$ by transform $\xi(u)$.
Recall the left-invariant vector fields $A_1,A_2,A_3$ out of which operators $\partial_+,\partial_-,\partial_0$
were defined as linear combinations.
Thus we are interested in commutators
$[\sigma_{A_j}(x,\xi),\sigma_A(x,\xi)]$,
and we already know the good behaviour of the case $A_3$, since $\partial_0={\rm i}A_3$.
Now,
\begin{eqnarray*}
    \sigma_{A_1}(x,\xi) & = & t^\xi(v_1)^\ast
    \ \sigma_{A_3}(x,\xi)\ t^\xi(v_1), \\
    \sigma_{A_2}(x,\xi) & = & t^\xi(v_2)^\ast\ \sigma_{A_3}(x,\xi)\ t^\xi(v_2), \end{eqnarray*} where $v_j\in G$ and $t^\xi:{\rm SU}(2)\to {\rm U}(2\xi+1)$ is the usual irreducible unitary matrix representation.
Hence
\begin{eqnarray*}
    \left[ \sigma_{A_1}(x,\xi),\sigma_A(x,\xi) \right]
    & = & \left[
    t^\xi(v_1)^\ast \sigma_{A_3}(x,\xi) t^\xi(v_1), \sigma_A(x,\xi) \right] \\
    & = & t^\xi(v_1)^\ast
    \left[ \sigma_{A_3}(x,\xi), t^\xi(v_1) \sigma_A(x,\xi) t^\xi(v_1)^\ast
    \right]
    t^\xi(v_1).
\end{eqnarray*}
We may assume that here $A\in\Psi^m(G)$ is left-invariant, i.e.
$A=(f\mapsto f\ast a)$, so that
$\sigma_A(x,\xi)=\widehat{a}(\xi)$ for every $\xi\in\Bbb N/2$ (in these symbol commutator properties, we may always fix $x\in G$, even if $A$ is not left-invariant).
Now $\left((x,\xi)\mapsto \widehat{a}(\xi)\right)\in S^m(G)$, and due to the diffeomorphism-invariance of $\Psi^m(G)$, also $$
    (x,\xi)\mapsto t^\xi(v_1)\ \widehat{a}(\xi)\ t^\xi(v_1)^\ast $$ must belong to $S^m(G)$, thus ``decaying rapidly off-diagonal''.
Thereby we already understand the behaviour of the commutator $$
    \left[ \sigma_{A_3}(x,\xi), t^\xi(v_1) \sigma_A(x,\xi) t^\xi(v_1)^\ast
    \right],
$$
and this is similar (by transform $t^\xi(v_1)^\ast$, see above) to $$
    \left[ \sigma_{A_1}(x,\xi),\sigma_A(x,\xi) \right].
$$
The same kind of reasoning applies to the commutator $$
    \left[ \sigma_{A_2}(x,\xi),\sigma_A(x,\xi) \right].
$$
Thus the symbol commutator condition is satisfied if we require that $$
    (x,\xi)\mapsto t^\xi(u)\ \sigma_A(x,\xi)\ t^\xi(u)^\ast $$ ``decays rapidly off-diagonal'' for every $u\in G$.

\section{Global calculus}

Properties (a$_{\Bbb R^n}$) and (b$_{\Bbb R^n}$) can be used to construct the global calculus of the introduced matrix symbols.
For example, we have the well-known composition formula for pseudo-differential operators $A$ and $B$ in $\Rn$:
$$
    \sigma_{AB}(x,\xi) \sim \sum_{\alpha\geq 0} \frac{1}{\alpha!}
           \ (\partial_\xi^\alpha\sigma_A)(x,\xi)
           \ (\partial_x^\alpha\sigma_B)(x,\xi),
$$
see e.g. \cite{KN65}, \cite{Ho85}, \cite{Ku81}.
On the torus ${\Bbb T^n}$, properties (a$_{\Bbb T^n}$) and (b$_{\Bbb T^n}$) allow to describe the composition of pseudo-differential operators in terms of toroidal symbols in an analogous way:
$$
    \sigma_{AB}(x,\xi) \sim \sum_{\alpha\geq 0} \frac{1}{\alpha!}
           \ (\triangle_\xi^\alpha\sigma_A)(x,\xi)
           \ (\partial_x^{(\alpha)}\sigma_B)(x,\xi).
$$
This formula was obtained in \cite{TV98} and was further extended to compositions with toroidal Fourier integral operators in \cite{RT07}.
In the case of ${\rm SU}(2)$, we can use properties (a$_G$) and (b$_G$) to show that it also has the form $$
    \sigma_{AB}(x,\xi) \sim \sum_{\alpha\geq 0} \frac{1}{\alpha!}
           \ (\triangle_\xi^\alpha\sigma_A)(x,\xi)
           \ (\partial_x^{(\alpha)}\sigma_B)(x,\xi).
$$
In fact, this formula can be also viewed in the following way. Let $$
     Q^\alpha\sigma_A(x)=\pi (y\mapsto \check q_\alpha(y)\ s_A(x)(y)) $$ with $\check q_\alpha\in C^\infty(G)$ satisfying $\check q_\alpha(\exp(z)) = z^{\alpha}$ for $z$ in a small neighbourhood of the origin on the Lie algebra $\mathfrak g$ of $G$.
Then we have the following elements of the calculus:
\begin{thm}[M. E. Taylor \cite{Ta84}]
\begin{eqnarray*}
    \sigma_{A^\ast}(x) & \sim & \sum_{\alpha\geq 0} \frac{1}{\alpha!}
          \ Q^\alpha\ \partial_x^\alpha\ \sigma_A(x)^\ast,\\
     \sigma_{AB}(x) & \sim & \sum_{\alpha\geq 0} \frac{1}{\alpha!}
          \ \left( Q^\alpha \sigma_A(x)\right)
          \ \partial_x^\alpha \sigma_B(x).
\end{eqnarray*}
\end{thm}

\begin{thm}
Let $A\in\Psi^m(G)$ with expansion
$A\sim\sum_{j=0}^\infty A_j$, $A_j\in\Psi^{m-j}(G)$, and assume that $x\mapsto\sigma_{A_0}(x)^{-1}$ is a symbol of $B_0\in\Psi^{-m+1-\varepsilon}(G)$ (for some $\varepsilon > 0$).
Then $A$ is elliptic with a parametrix $B$, $\sigma_B\sim\sum_{k=0}^\infty\sigma_{B_k}$, where \begin{eqnarray*}
     \sigma_{B_0}(x) & = & \sigma_{A_0}(x)^{-1},\\
     \sigma_{B_N}(x) & = & -\sigma_{B_0}(x)\sum_{k=0}^{N-1}
     \sum_{j=0}^{N-k}
           \sum_{\gamma:j+k+|\gamma|=N}
           (Q^\gamma \sigma_{A_j}(x))
           \ \partial_x^{\gamma}\sigma_{B_k}(x).
\end{eqnarray*}
\end{thm}

For further details of these results we refer to \cite{Ta84}, \cite{Tu01}, \cite{Tu04}, and for the relation to introduced symbol classes we refer to \cite{RT08}.



\end{document}